\documentclass[conference]{IEEEtran}
\IEEEoverridecommandlockouts

\usepackage{cite}
\usepackage{amsmath,amssymb,amsfonts,mathtools}
\usepackage{algorithmic}
\usepackage{graphicx}
\usepackage{textcomp}
\usepackage{xcolor}
\def\BibTeX{{\rm B\kern-.05em{\sc i\kern-.025em b}\kern-.08em
    T\kern-.1667em\lower.7ex\hbox{E}\kern-.125emX}}
\begin{document}

\title{Sequential Change Detection under Markov Setup with Unknown Pre-change and Post-change Distributions} 

\author{\IEEEauthorblockN{A. B. Gulaguli}
\IEEEauthorblockA{\textit{Dept Of Electrical Engineering} \\
\textit{Shiv Nadar University}\\
Delhi NCR, India \\
ag497@snu.edu.in}
\and
\IEEEauthorblockN{S. Singh}
\IEEEauthorblockA{\textit{Dept Of Electrical Engineering} \\
\textit{Shiv Nadar University}\\
Delhi NCR, India \\
ss632@snu.edu.in}
\and
\IEEEauthorblockN{R. K. Bansal}
\IEEEauthorblockA{\textit{Dept Of Electrical Engineering} \\
\textit{Shiv Nadar University}\\
Delhi NCR, India \\
rakesh.bansal@snu.edu.in}
}
\maketitle
\begin{abstract}
In this work we extend the results developed in 2022 for a sequential change detection algorithm making use of Page’s CUSUM statistic, empirical distribution as an estimate of pre-change distribution and a universal code as a tool for estimating the post-change distribution, from i.i.d. to Markov setup.
\end{abstract}

\section{Introduction}
Change detection has numerous applications including industrial quality control, and fault detection in computer communication networks and has been extensively studied in literature. Sequential hypothesis testing and sequential change detection have a long history.

Wald formalized the sequential hypothesis testing problem in 1945 \cite{wald1945}. He proposed the Sequential Probability Ratio Test (SPRT) and proved its optimality along with Wolfowitz \cite{waldwolfowitz1948}. Lai proved asymptotic optimality of invariant SPRT for classes of sources that satisfy a stability condition for likelihood ratio \cite{lai1981}.

Page introduced the Cumulative Sum Test (CUSUM) to address the problem of sequential change detection in 1954 \cite{page1954}. Lorden proved the asymptotic optimality of CUSUM for memoryless sources in minimax framework \cite{lorden1971} while Moustakides proved its exact optimality \cite{moustakides1986}. Shiryaev \cite{shiryaev1963,shiryaev1978} and Pollak \cite{Pollak:1985} proposed distinct approaches for the change detection problem under Bayesian and non-Bayesian framework respectively. Lorden's result was extended for sources with memory by Bansal and Kazakos \cite{bansal1986} who assume independence of observations before and after change which was further weakened in Lai's extension \cite{lai1998}. \cite{willsky1976,laishan1999} are other works addressing the unknown post-change parameter case.

Page's CUSUM Test requires that both  pre-change and post-change distributions should be known which is not always the case in real world situations. Jacob and Bansal \cite{jacobbansal2008} used ideas from information theory to develop a modified SPRT test as well as modified CUSUM test for the case where post-change distribution is unknown and a universal code is used in the place of log-likelihood of post-change distribution for memoryless sources. They showed the optimality of the modified test in asymptotic regime where the rate of false alarm goes to 0. They extended the idea of modified test for sources with memory under Berk's mixing conditions \cite{jacob2010}.

Malik and Bansal \cite{malikbansal2022} considered the case where both pre-change and post-change distributions are unknown for memoryless sources. They further modified the test introduced by Jacob and Bansal \cite{jacobbansal2008}, by using an empirical estimate in the place of true pre-change distribution without affecting the post-change distribution estimate and proved the asymptotic optimality of the proposed test.

In this paper, we extend the modification demonstrated by Malik and Bansal \cite{malikbansal2022} for the more general class of stationary and ergodic finite-order Markov sources that satisfy Berk's mixing conditions\cite{1973} and make use of  Bansal and Kazakos' extension \cite{bansal1986} to study the performance of proposed test and establish its optimality in the asymptotic regime where the rate of false alarm goes to 0. We evaluate the performance of the proposed test by providing an upper bound on its probability of error, showing that it will terminate with probability one under the alternative hypothesis, and studying a cost function for the number of samples required to achieve a specific probability of error as well as a cost function for the detection delay required to achieve a specific false alarm rate.

Section II introduces the preliminaries required for the rest of the paper. In Section III, we introduce the empirical estimate, the modified test and study the behaviour of proposed test under pre-change and post-change distribution. In Section IV, we study the performance of modified test and prove its asymptotic optimality and conclude the paper in Section V.
\section{Preliminaries}
Let \( X = \{X_n\}_{n=1}^\infty \) be a random process (source) taking values from a finite set called the source alphabet $\mathcal{X}$. A sequence of source symbols \( x_i, x_{i+1}, \dots, x_j \) is denoted by \( x_i^j \) and corresponding random variables  \( X_i, X_{i+1}, \dots, X_j \) is denoted by \( X_i^j \). All the logarithms will be taken to the base 2 in the remaining, unless specified otherwise.

For a fixed \textit{n}, a fixed-to-variable code (FV code) is a one-to-one mapping \(\phi_n : X^n \to \{0, 1\}^*\), where \(\{0, 1\}^*\) denotes the set of all finite-length strings. Let \( L(x_1^n) = |\phi_n(x_1^n)|\) denote the length of the codeword for a given sequence \( x_1^n \) which satisfies prefix condition. Therefore, Kraft inequality holds, given as,
\[ \sum_{x_1^n \in X^n} 2^{-L(x_1^n)} \leq 1. \tag{1}\]

Let \( \mu \) be the distribution of \( X = \{X_n\}_{n=1}^\infty \). The entropy rate for stationary \( \mu \)  is well defined \cite{cover2006}. For a stationary ergodic process, entropy rate of the source is a well-known lower bound for the performance of any source coding algorithm \cite{cover2006}. Optimal codes reach this lower bound. A FV code with length function \( L(x_1^n) \) is said to be \textit{pointwise universal} if,
\[\frac{1}{n} \left( L(x_1^n) + \log \mu(x_1^n) \right) \to 0 \quad \mu\text{ a.s.} \tag{2}
\]  
For stationary ergodic sources, the existence of such codes is guaranteed \cite{csiszar2004,cover2006}. For a class of sources \( \mathcal{M} \), a FV code with length function \( L(x_1^n) \) is said to be \textit{strongly pointwise universal} if:
\[R_L^n = \sup_{\mu \in \mathcal{M}} \max_{x_1^n \in X^n} \left( L(x_1^n) + \log \mu(x_1^n) \right) \sim o(n). \tag{3}\]
The condition for pointwise universality (2) is equivalent to having \( \mu (J_\epsilon < \infty) = 1 \) for all \( \epsilon > 0 \), where,
\begin{equation}
J_\epsilon = \sup_{n \geq 1} 
\left\{ 
n : 
\left| 
\frac{1}{n} \left( L(x_1^n) + \log \mu(x_1^n)\right)
\right| 
\geq \epsilon 
\right\}.
\tag{4}
\end{equation}

A stronger version of the convergence in condition (2) requires that,
\begin{equation}
\mathbb{E}_\mu (J_\epsilon) < \infty 
\quad \text{for all} \quad \epsilon > 0.
\tag{5}
\end{equation}
Strongly pointwise universal codes (3) satisfy this condition which is called 1-quick universal by Jacob \cite{jacob2010}, for every source in the class \( \mathcal{M} \). For Markov sources upto a fixed order, the existence of codes that satisfy conditions (3) and (5) are known (Lempel-Ziv codes (LZ-77, LZ-78) \cite{zivlempel1977,zivlempel1978} and the Context Tree Weighting code \cite{willems1995} )  and we consider such codes only in our paper.
The divergence rate for two stationary processes $\mu$ and $\nu$ is defined as,
\[ D(\mu || \nu) = \lim_{n \to \infty} \frac{1}{n} \sum_{x_1^n \in X^n} \mu(x_1^n) \log \frac{\mu(x_1^n) }{\nu(x_1^n) }\tag{6} \]
 provided the limit exists.
 
\textbf{Theorem 1 (Barron)\cite{barron1985,gray1990}}: Let \( \mu \) be a stationary, ergodic process and \( \nu \) be a stationary, finite-order Markov process such that \( \mu \ll \nu \). Then,
\[\frac{1}{n} \log \frac{\mu(x_1^n)}{\nu(x_1^n)} \to D(\mu || \nu) \quad \mu \text{ a.s.} \quad \tag{7} \]
where \( D(\mu || \nu) \), the divergence rate is well defined.

The pointwise convergence in the above theorem is equivalent to having $\mu (K_\epsilon < \infty) = 1 \quad \text{for all} \quad \epsilon > 0$ where,
\[K_\epsilon = \sup_{n \geq 1} \left\{ n : \left| \frac{1}{n} \log \frac{\mu(x_1^n)}{\nu(x_1^n)} - D(\mu || \nu) \right| \geq \epsilon \right\}\quad. \tag{8} \]
A stronger version of this convergence, called 1-quick convergence of log-likelihood ratios according to Lai's terminology \cite{lai1981,jacob2010}, requires that,
\[\mathbb{E}_\mu (K_\epsilon) < \infty \quad \text{for all} \quad \epsilon > 0. \quad \tag{9} \]

For a fixed sample hypothesis test with the null hypothesis $H_0 : \mu = \mu_0$  and the alternative $H_1 : \mu = \mu_1$, consider the discriminant function,
\[
  h(x_1^n) = \log \mu_1(x_1^n) - \log \mu_0(x_1^n) - n\lambda. \tag{10}
\]
After observing $n$ samples $x_1^n$ we reject the null hypothesis if $h(x_1^n) \geq 0$. The performance of this test is characterized by the following lemma.

\textbf{Lemma 1 \cite{chernoff1952,ziv1988}:} If $\mu_0$ and $\mu_1$ are stationary, finite-order Markov process and stationary, ergodic process such that $\mu_1 \ll \mu_0$ and $D(\mu_1 \| \mu_0) > \lambda$. For the test induced by the discriminant function $h(x_1^n)$ we have:
\begin{align}
    \mu_0\big(h(x_1^n) \geq 0 \big) &\;\leq\; 2^{-n\lambda} \tag{11} \\
    \lim_{n \to \infty} \mu_1\big(h(x_1^n) \geq 0 \big) &= 1.
    \tag{12}
\end{align}
Equation (12) remains true even if the threshold $0$ is replaced by any other constant.

\textbf{Berk's Mixing Conditions \cite{1973,bansal1986,jacob2010}:} Let
\[
p_n(\Delta) = \mu_1 \left( \frac{1}{n} \log \frac{\mu_1(x_1^n)}{\mu_0(x_1^n)} < \Delta \right) \tag{13}
\]

and

\[
q_n(\Delta) = \mu_0 \left( \frac{1}{n} \log \frac{\mu_0(x_1^n)}{\mu_1(x_1^n)} < \Delta \right). \tag{14}
\]

For $\Delta < D(\mu_1 \| \mu_0)$,
\[
\lim_{n \to \infty} n p_n(\Delta) = 0 
\quad \text{and} \quad 
\sum_{n \geq 1} p_n(\Delta) < \infty. \tag{15}
\]

For $\Delta < D(\mu_0 \| \mu_1)$,
\[
\lim_{n \to \infty} n q_n(\Delta) = 0 
\quad \text{and} \quad 
\sum_{n \geq 1} q_n(\Delta) < \infty. \tag{16}
\]

\subsection{Previous Work}
For the case where both pre-change distribution $\mu_0$ and post-change $\mu_1$ distribution are known, Page proposed a test \cite{page1954} which is given as follows.  

\textbf{Page’s CUSUM Test (Test 1) \cite{page1954}:} Starting with $n=1$, for samples $x_1, x_2, \ldots, x_n$, and $\gamma > 1$, if  
\begin{equation}
\max_{k \leq n} \left( \log \mu_1(x_k^n) - \log \mu_0(x_k^n) \right) \geq \log \gamma,
\tag{17}
\end{equation}
stop the test and decide that a change has occurred. Else, sample for $x_{n+1}$ and continue the test.  

\medskip

\textbf{Jacob–Bansal Page CUSUM Test (Test 2)\cite{jacobbansal2008,jacob2010}:} A modification of Page's CUSUM test was proposed for cases where the post-change distribution $\mu_1$ is unknown for i.i.d. sources and later extended for sources with memory satisfying (9), (15), (16). The modified CUSUM Test, where a strongly pointwise universal code is used to estimate the post-change distribution, is given as follows.  

Starting with $n=1$, for samples $x_1, x_2, \ldots, x_n$, $\gamma > 1$ and $\lambda < D(\mu_1 \Vert \mu_0)$, if  
\begin{equation}
\max_{k \leq n} \left( -L(x_k^n) - \log \mu_0(x_k^n) - n \lambda \right) \geq \log \gamma,
\tag{18}
\end{equation}
stop the test and decide that a change has occurred. Else, sample for $x_{n+1}$ and continue the test.  

Authors introduced $\lambda$ here to ensure a negative drift in the test statistic in the pre-change regime, as is the case in the classical Page’s test in the parametric setup.  

\medskip

An auxiliary stopping time for this test is defined as
\begin{equation}
N_0(\alpha) = \inf \left\{ n : -L(x_1^n) - \log \mu_0(x_1^n) - n\lambda \geq -\log \alpha \right\}. \tag{19}
\end{equation}

An error occurs when Test 2 stops even though no change in distribution has occurred. The probability of error for the test is calculated using the stopping time $N_0(\alpha)$.

\textbf{Malik–Jacob-Bansal Page CUSUM Test (Test 3) \cite{malikbansal2022}:} Since JB-Page’s CUSUM test requires the knowledge of pre-change distribution, a modification of this test was proposed for cases where the pre-change distribution $\mu_0$ is unknown for i.i.d. sources. The modified CUSUM test, where a training sequence is used to estimate the pre-change distribution, is given as follows.

Test (\(\gamma>1\), \(\lambda>0\)) : Start with $n=1$, for samples \(x_{n_0+1}, x_{n_0+2}, \dots, x_{n_0+n}\). If: 
\begin{equation}
\max_{1 \leq k \leq n} \left( - L(x_{n_0+k}^{n_0+n}) - \log \hat{\mu}_0({x}_{n_0+k}^{n_0+n}) - n\lambda \right) \geq \log \gamma,
\tag{20}
\end{equation}
stop the test and decide that a change has occurred. Otherwise, sample for \(x_{n_0+n+1}\) and continue the test.

Readers can refer to \cite{malikbansal2022} to understand why universal code cannot be used to estimate pre-change distribution. 

\section{Modified CUSUM Test For Markov}
We consider stationary, ergodic, mutually independent, first-order
Markov sources $\mu_0$ and $\mu_1$ satisfying (9), (15), (16) and $\mu_1\ll\mu_0$  as well as universal codes satisfying (5) for analysis. We also assume that if $\mu_0\in\mathcal{M}_0$ and $\mu_1\in\mathcal{M}_1$ then $\min_{\substack{\mu_0 \in \mathcal{M}_0 \\ \mu_1 \in \mathcal{M}_1}} 
D(\mu_1 \Vert \mu_0) = \lambda_0 > 0$.
\subsection{Estimating Pre-change Distribution}
For a stream of source symbols with no change up to \( n_0 \), the test has two stages. In the first stage, we estimate the pre-change distribution using the first \( n_0 \) symbols, written as \( x_1, x_2, \dots, x_{n_0} \) and the subsequent symbols as \( x_{n_0+1}, x_{n_0+2}, \dots, x_{n_0+n} \). The prechange distribution estimate, \( \hat{\mu}_0 \), is calculated from initial sample of \( n_0 \) symbols as,
   \[\hat{\mu}_0(b,a) = \frac{N(b,a) \mid x_1^{n_0}}{n_0 - 1} \quad \forall (b,a) \in \mathcal{X}^2\tag{21}\] 
   \[\hat{\mu}_0(b) = \frac{N(b) \mid x_1^{n_0}}{n_0} \quad \forall b \in \mathcal{X} \tag{22}\]
   \[\hat{\mu}_0(x_1^{n_0}) = \frac{\hat{\mu}_0(x_1, x_2) \hat{\mu}_0(x_2, x_3) \dots \hat{\mu}_0(x_{n_0-1}, x_{n_0})}{\hat{\mu}_0(x_2) \dots \hat{\mu}_0(x_{n_0-1})} \tag{23}\]
   \[\hat{\mu}_0(x_1^{n_0}) = \hat{\mu}_0(x_1) \prod_{i=2}^{n_0} \hat{\mu}_0(x_i / x_{i-1}) \tag{24} \]  
 
Here, \( N(b) \mid x_1^{n_0} \) is the number of times symbol \( b \) appears in the sequence \( x_1^{n_0} \), whereas \( N(b,a) \mid x_1^{n_0} \) is the number of times the symbols  \((b,a)\) jointly appear in the sequence \( x_1^{n_0} \). Hence,
\( \hat{\mu}_0(b) \) is the fraction of times any symbol \(b\in \mathcal{X}\) appears in the sequence \( x_1^{n_0} \), whereas \( \hat{\mu}_0(b,a) \) is the fraction of times any symbol pair\( (b,a)\in \mathcal{X}^2\) appear in the sequence \( x_1^{n_0} \) jointly.

We assume that \( n_0\) is large enough so that each symbol $b \in \mathcal{X}$ which has non-zero probability under ${\mu}_0$ has appeared at least once in \( x_1^{n_0} \) and each symbol pair $(b,a) \in \mathcal{X}^2$ which has non-zero probability under ${\mu}_0$ has appeared at least once in \( x_1^{n_0} \) (else $\hat{\mu}_0$ is not an accurate estimate of ${\mu}_0$).

Since ${\mu}_0$ is a stationary and ergodic first-order Markov source, by Ergodic Theorem, we can say that there exists an \( n_0\) such that,
\[
\mu_0 \left( \left\{ x_1^{n_0} : \big| \hat{\mu}_0(a \mid b) - \mu_0(a \mid b) \big| > \delta \right\} \right) \leq \epsilon, \tag{25}
\]
for one or more \((b,a) \in \mathcal{X}^2.\)

\[\mu_0 \left( \{ x_1^{n_0} : | \hat{\mu}_0(b) - \mu_0(b) | > \delta' \} \right) \leq \epsilon', \tag{26}\] for one or more $b \in \mathcal{X}$

\subsection{Test}
Using the estimate \(\hat{\mu}_0\), we extend the modified CUSUM Test for Markov sources by defining it as,

Test ($\gamma>1$): Start with $n=1$ for samples \(x_{n_0+1}, x_{n_0+2}, \dots, x_{n_0+n}\). If
\[\max_{1 \leq k \leq n} \left( - L(x_{n_0+k}^{n_0+n}) - \log \hat{\mu}_0({x}_{n_0+k}^{n_0+n}) - n\lambda \right) \geq \log \gamma, \tag{27}\]
stop the test and decide that a change has occurred. Else, sample for \(x_{n_0+n+1}\) and continue the test. The upward drift rate and downward drift rate of this test are \(D(\mu_1 || \hat{\mu}_0) - \lambda\) and \(\lambda - D(\mu_0 || \hat{\mu}_0)\) respectively. Hence, for this test to work, we need:
\[D(\mu_0 || \hat{\mu}_0) < \lambda < D(\mu_1 || \hat{\mu}_0) \tag{28}\]
which implies,
\[D(\mu_0 || \hat{\mu}_0) < D(\mu_1 || \hat{\mu}_0). \tag{29}\]
Auxiliary stopping time for this modified test is defined as,
\begin{equation}
\begin{aligned}
N(\alpha)
= \inf \Bigl\{ n :
\Bigl(
&- L(x_{n_0+1}^{n_0+n})
- \log \hat{\mu}_0(x_{n_0+1}^{n_0+n}) \\
&- n\lambda
\Bigr)
\ge -\log \alpha
\Bigr\}.
\end{aligned}
\tag{30}
\end{equation}
\subsection{Behaviour Under Pre-change Distribution}
Introducing the following notation to study the expected behaviour of the auxillary stopping time $N(\alpha)$ under pre-change distribution \( \mu_0 \) and post-change distribution \( \mu_1 \).
\[f_n(x_{n_0+1}^{n_0+n}) \triangleq \frac{\log \frac{\mu_0(x_{n_0+1}^{n_0+n})}{\hat{\mu}_0(x_{n_0+1}^{n_0+n})}}{n}
\]

\[l_n(x_1^{n_0+n}) \triangleq -L(x_{n_0+1}^{n_0+n}) - \log\hat{\mu}_0(x_{n_0+1}^{n_0+n}) - n\lambda
\]
\[
\begin{aligned}
B_n \triangleq \{x_1^{n_0+n} : l_k(x_1^{n_0+k}) < -\log \alpha, \; 1 \leq k < n, \; l_n(x_1^{n_0+n}) \\ \geq -\log \alpha\}     
\end{aligned}  
\]
\[
B \triangleq \bigcup_{n=1}^\infty B_n
\]
\[
\begin{aligned}
E_n \triangleq \Bigl\{ x_1^{n_0+n} : 
&\lvert\hat{\mu}_0(a \mid b) - \mu_0(a \mid b) \big| < \delta \quad \forall \, (a,b) \in \mathcal{X}^2, \\
&\lvert\hat{\mu}_0(b) - \mu_0(b) \big| < \delta' \quad \forall \, b \in \mathcal{X} 
\Bigr \}
\end{aligned}
\]
\[
\text{ where } \hat{\mu}_0 \text{ is based on samples } x_1^{n_0}
\]

\[E \triangleq \bigcup_{n=1}^\infty E_n\]

\textbf{Lemma 2:} If \( x_1^{n_0+n} \in E_n \), then \( f_n(x_{n_0+1}^{n_0+n}) \) is bounded above for all \( n \) as  
$f_n \leq W_n=\frac{(n-1)\log (\beta_1)+\log (\beta_2)}{n}$, where $\beta_1=\frac{\mu_0(a' \mid b')}{\mu_0(a' \mid b') - \delta}, \beta_2=\frac{\mu_0(b')}{\mu_0(b') - \delta'}$ and \((b',a')\) is such that \((b',a')\) is the minimum transition probability symbol pair and \(b'\) is such that \(b'\) is the minimum probability symbol. As \( \delta \to 0 \), \( \delta' \to 0 \) and \( n_0 \to \infty \), we have \( \beta_1 \to 1 \) and \( \beta_2 \to 1 \).

\textit{Proof:} As $x_i$'s follow first order Markov process, we can write $f_n\big(x_{n_0+1}^{n_0+n}\big)$ as,
\[
f_n \big(x_{n_0+1}^{n_0+n} \big) 
= \sum_{(b,a) \in \mathcal{X}^2} 
\frac{N(b,a \mid x_{n_0+1}^{n_0+n})}{n} 
\cdot \log \frac{\mu_0(a \mid b)}{\widehat{\mu}_0(a \mid b)}
\]
\[
+\sum_{b \in \mathcal{X}} 
\frac{N(b \mid x_{n_0+1})}{n} 
\cdot \log \frac{\mu_0(b)}{\widehat{\mu}_0(b)}. 
\tag{31}
\]
As $x_{n_0+1}^{n_0+n} \in E_n$, for \(x_i\)'s sitting at any positions for \(i = n_0+1, \ldots, n_0+n\),
\[\frac{\mu_0(a\mid b)}{{\widehat{\ \mu}}_0(a\mid b)} <\frac{\mu_0(a\mid b)}{\mu_0(a\mid b)-\delta},\frac{\mu_0(b)}{{\widehat{\ \mu}}_0(b)} <\frac{\mu_0(b)}{\mu_0(b)-\delta'}.\tag{32} \] 
The terms on RHS are decreasing functions of $\mu_0(a \mid b)$ and $\mu_0(b)$ respectively, So if $b'$ is the minimum probability symbol and $(b',a')$ is the minimum transition probability symbol pair then we can bound the above expression for \(f_n (x_{n_0+1}^{n_0+n})\) as (for $\delta\leq{\mu_0(a \mid b)}$ and $\delta'\leq{\mu_0(b)}$):
\[
\begin{aligned}
 f_n (x_{n_0+1}^{n_0+n} ) \leq \frac{n-1}{n} \cdot \log\left(\frac{\mu_0(a'\mid b')}{\mu_0(a'\mid b')-\delta}\right)+ \\ \frac{1}{n} \cdot\log(\frac{\mu_0(b')}{\mu_0(b')-\delta'}). 
\end{aligned}
\tag{33} \] 
\[f_n (x_{n_0+1}^{n_0+n} ) \leq \frac{(n-1)\log (\beta_1)+\log (\beta_2)}{n}=W_n. \tag{34} \]

 \textbf{Lemma 3:} The probability of error for the proposed test is bounded above by,
\[\mu_0(N(\alpha) < \infty) \le 
\frac{\alpha}{\dfrac{\beta_1 \left( 2^{\lambda - \log \beta_1} - 1 \right)}{\beta_2}} 
+ \epsilon_0 \tag{35}\]
where \(\epsilon_0 \to 0\), \(\beta_1 \to 1\) and \(\beta_2 \to 1\) as \(n_0 \to \infty\). 

\textit{Proof:} Using the definition of $f_n(x_{n_0+1}^{n_0+n})$,
\[
\mu_0\!\left(x_{n_0+1}^{n_0+n}\right) 
= \hat{\mu}_0\!\left(x_{n_0+1}^{n_0+n}\right) \cdot 2^{\,n f_n(x_{n_0+1}^{n_0+n})}. \tag{36}
\]

The probability of error is given as,
\[
\begin{aligned}
\mu_0\!\left(N(\alpha) < \infty\right) 
= \mu_0(B)=\sum_{n=1}^{\infty} \mu_0 \!\left( B_n \cap E_n \right) 
   + \\ \mu_0 \!\left( B \cap (E)^c \right).   
\end{aligned}
\tag{37}
\]
\[
= \sum_{n=1}^{\infty} 
   \left( \sum_{x_{1}^{n_0+n} \in B_n \cap E_n}
   \mu_0 \!\left( x_{1}^{n_0+n} \right) \right) 
   + \mu_0 \!\left( B \cap (E)^c \right). \tag{38}
\]

As $B_n$ consists of only $n$-long sequences, we have $B_n \cap E = B_n \cap E_n $. For first term of (38) we have,
\[
\sum_{x_{1}^{n_0+n} \in B_n \cap E_n}
\mu_0 \bigl( x_{1}^{n_0+n} \bigr) \leq \sum_{x_{1}^{n_0+n} \in B_n \cap E_n} 
\mu_0 \bigl( x_{n_0+1}^{n_0+n} \bigr) \tag{39}
\]
\[
\leq \sum_{x_{1}^{n_0+n} \in B_n \cap E_n} 
(\hat{\mu}_0\!\left(x_{n_0+1}^{n_0+n}\right) \cdot 2^{\,n f_n(x_{n_0+1}^{n_0+n})}\bigr) \tag{40}
\]
\[
\leq \sum_{x_{1}^{n_0+n} \in B_n \cap E_n} 
(\hat{\mu}_0\!\left(x_{n_0+1}^{n_0+n}\right) \cdot 2^{\,n W_n}\bigr) \tag{41}
\]  
\[
\leq \sum_{x_{1}^{n_0+n} \in B_n \cap E_n} 2^{\log \alpha-L(x_{n_0+1}^{n_0+n})\,-n(\lambda-W_n) } \tag{42}
\]
\[
=\alpha \cdot 2^{\,-n(\lambda-W_n)}(\sum_{x_{1}^{n_0+n} \in B_n \cap E_n} 2^{-L(x_{n_0+1}^{n_0+n})}) \tag{43}
\]
\[
\leq\alpha \cdot 2^{\,-n(\lambda-W_n)}. \tag{44}
\]

(40) follows from (36) while the third inequality is due to Lemma  2 as \( x_1^{n_0+n} \in E_n \). The fourth inequality comes from definition of $B_n$ while the last inequality comes from Kraft's inequality.
As $P(A \cap B) \leq{P(B)}$, and from (25) and (26), we have,
\[
\mu_0(B \cap (E)^c)\leq{\mu_0((E)^c)}\leq{\epsilon}+{\epsilon'}={\epsilon_0}.
\tag{45}\]

Combining both equations and from (34) we have,
\[
\mu_0(N(\alpha) < \infty) \le {\alpha}\sum_{n=1}^{\infty} 2^{-n (\lambda - \log\beta_1)+\log(\frac{\beta_2}{\beta_1})} + \epsilon_0 \tag{46}
\]
\[
\mu_0(N(\alpha) < \infty) \le 
\frac{\alpha}{\dfrac{\beta_1 \left( 2^{\lambda - \log \beta_1} - 1 \right)}{\beta_2}} 
+ \epsilon_0.
\tag{47}
\]

By comparing the probability of error expressions in Lemma 3 and for known pre-change case (18), we observe that the probability of error increases for our test. The probability of error now becomes a function of \(\log \beta_1\), \(\beta_1\) and \(\beta_2\).
For the probability of error to be bounded, we require:  
\[\lambda > \max(\log \beta_1, D(\mu_0 || {\widehat{\mu}}_0)). \tag{48}\]  
In contrast, for known pre-change case (18), \(\lambda > 0\) would have sufficed. However, as \(n_0 \to \infty\), the lower bound \(\max(\log \beta_1, D(\mu_0 || {\widehat{\mu}}_0))\), the error term \(\epsilon_0\) and \(\beta_1,\beta_2\) can be made as small as needed.  

Hence, \[\max (\log \beta_1, D(\mu_0 || {\widehat{\mu}}_0)) < \lambda < D(\mu_1 || \hat{\mu}_0). \tag{49}\]

\subsection{Behaviour Under Post-change Distribution}
\textbf{Lemma 4:} Let the conditions of Lemma 1 hold true. The proposed test (27) will terminate with probability 1 under \(\mu_1\) if \(\lambda < D(\mu_1 || {\widehat{\mu}}_0)\), i.e.,
\[\mu_1(N(\alpha) < \infty) = 1. \tag{50}\]  

\textit{Proof:} The proof of this lemma is the same as that of Lemma 3 in \cite{jacobbansal2008}, since its proof is dependent on Lemma 1, replacing $\mu_0$ with ${\widehat{\mu}}_0$. The replacement is valid in this case as the behaviour of the test under $\mu_1$ is identical to that of a test where pre-change distribution is $\mu_0$.

\textbf{Theorem 2:} Let \(\mu_0\) and \(\mu_1\) be stationary, ergodic first-order Markov sources satisfying (9), (15), (16) and $\mu_1 \ll \mu_0$. For the Modified CUSUM Test (27) with a universal code satisfying (5) and stopping time \(N(\alpha)\) (30) as \(\alpha \to 0\), for \(\lambda < D(\mu_1 || {\widehat{\mu}}_0)\), we have, 
\[{E_1}(N(\alpha) \, | \, x_1^{n_0}) \sim \frac{| \log \alpha |}{D(\mu_1 || {\widehat{\mu}}_0) - \lambda}. \tag{51}\]  

\textit{Proof:} This theorem can be proved by arguing along similar lines as Theorem 2 of \cite{jacobbansal2008} since its proof is dependent on (3), (9) which are satisfied by the codes and sources considered in this paper.

\section{Performance Bounds}
To study the performance of the proposed modified test, we use $M(\gamma)$ to denote the stopping time of the test,
\begin{equation}
\begin{aligned}
M(\gamma)
= \inf \Bigl\{ n : \max_{1 \le k \le n}
\Bigl(
&- L(x_{n_0+k}^{n_0+n})
- \log \hat{\mu}_0(x_{n_0+k}^{n_0+n}) \\
&- n\lambda
\Bigr)
\ge \log \gamma
\Bigr\}.
\end{aligned}
\tag{52}
\end{equation}
We work with conditional expectation while studying performance under \(\mu_1\), conditioned on the initial \(n_0\)-length segment generated by \(\mu_0\) as we used the empirical distribution \(\hat{\mu}_0\) generated by this segment in place of \(\mu_0\).

Consider the case where first \(m-1\) samples are generated by \(\mu_0\) and the rest of the samples come from \(\mu_1\). In this case, the change point is \(m\). Let \({P_m}\) and \({E_m}\) be the corresponding probability distribution and expectation for this case. Then, $P_m\big(x_{1}^{\,n_0+n}\big) = \mu_0\big(x_{1}^{\,m-1}\big) \, \mu_1\big(x_{m}^{\,n_0+n}\big)$ due to the mutual independence of sources before and after the change \cite{bansal1986,jacob2010}.

For a given stopping time $S$, Lorden proposed the following quantity (53) to serve as  minimax type of criterion for quickness of reaction to a change \cite{lorden1971}.
\[\overline{E_1}( S ) = \sup_{m \geq 1} \, \text{ess sup} \, {E_m} \left[ (S - m + 1)^+ \left| X_1 , \dots , X_{m-1} \right| \right]. \tag{53}\]
To investigate the behaviour of \(\overline{E}_{1}(M(\gamma) \, | \, x_1^{n_0})\), we first study the characteristics of \({E_0}(M(\gamma))\), which is the expected time after which the algorithm detects a change and stops the test for the no source change case. Ideally, this has to be \(\infty\) as the algorithm should not detect a change when there isn’t any.

\textbf{Theorem 3 \cite{bansal1986}:} Lorden’s Theorem was extended for stationary and ergodic processes \( (\mu_0, \mu_1) \) satisfying Berk's mixing conditions (15), (16) such that,
\[{\mu_0}(N < \infty) \leq \alpha. \tag{54}\]
For \( k = 1, 2, \dots \), let \( N_k \) be the stopping variable obtained by applying \( N \) to \( X_k, X_{k+1}, \dots \), such that, 
\[N^\ast = \min \{ N_k + k - 1 \mid k = 1, 2, \dots \tag{55} \}.\]
Then \( N^\ast \) is an extended stopping variable with:
\[{E_0}(N^\ast) \geq \frac{1}{2} \left( 1 + \frac{1}{\alpha} \right) > \frac{1}{2\alpha} \tag{56}\]
and,
\[\overline{E_1}(N^\ast) \leq {E_1}(N). \tag{57}\]

Define \( S_\gamma \) to be the class of all sequential change detection tests' stopping times \( S \) satisfying:
\[{E_0}(S) \geq \frac{1+\gamma}{2} > \frac{\gamma}{2}. \tag{58}\]

\textbf{Theorem 4 \cite{bansal1986}}: For two stationary, ergodic and mutually independent stochastic processes with memory \( \mu_0, \mu_1 \) satisfying Berk's mixing conditions (15), (16), as \( \gamma \to \infty \), we have,

\[\inf_{S \in S_\gamma} \overline{E_1}(S) \sim \frac{\log \gamma}{D(\mu_1 || \mu_0)}. \tag{59}\]

\textbf{Theorem 5}: Let \(\mu_0\) and \(\mu_1\) be stationary, ergodic and mutually independent first-order Markov sources satisfying (9), (15), (16), $\mu_1 \ll \mu_0$ and \( \lambda < D(\mu_1 || \widehat{\mu_0})\). For the proposed test (27) using universal code satisfying (5) with stopping time $M(\gamma)$, we have,
\begin{equation}
\begin{aligned}
{E_0}\left(M(\gamma)\right) \geq \frac{1}{2} + \frac{\gamma}{\frac{2 \beta_2}{\beta_1(2^{\lambda - \log \beta_1} - 1)} + 2 \epsilon_0 \gamma}\\> \frac{\gamma}{\frac{2 \beta_2}{\beta_1(2^{\lambda - \log \beta_1} - 1)} + 2 \epsilon_0 \gamma},   
\end{aligned}
\tag{60}    
\end{equation}
and as \( \gamma \to \infty \), we have,
\[\overline{E_1}(M(\gamma) \mid x_1^{n_0}) \sim \frac{\log \gamma}{D(\mu_1 || \widehat{\mu_0}) - \lambda}. \tag{61}\]

\textit{Proof}: Equation (60) follows from (35) of Lemma 3 and (56) of Theorem 3 (restatement of Theorem 1 from \cite{bansal1986}).
Equation (61) follows from (57) of Theorem 3 (restatement of Theorem 1 from \cite{bansal1986}) and Theorem 2.

\textbf{Theorem 6:} Let \(\mu_0\) and \(\mu_1\) be stationary, ergodic and mutually independent first-order Markov sources satisfying (9), (15), (16) and $\mu_1 \ll \mu_0$. For the Modified CUSUM Test (27) with universal code satisfying (5), we can choose $\lambda$ and $\eta$ such that stopping time $M(\eta)=M\in S_\gamma$ and as \( \gamma \to \infty \) we have,  

\[
\overline{E_1}(M\mid x_1,x_2,x_3,...) \sim (1 + \vartheta) \inf_{S \in S_\gamma} \overline{E_1}(S) \quad \text{a.s.} \ [\mu_0], \tag{62}
\]
where \( \vartheta \) can be chosen arbitrarily close to \( 0 \).  

\textit{Proof:} Let 
\[
\eta = \gamma\left( \frac{\beta_2}{\beta_1(2^{\lambda - \log(\beta_1)}-1)} + \epsilon_0 \gamma\right). \tag{63}
\]  

Now, if \( \epsilon_0 = \frac{1}{\gamma} \), then \( \eta \) can be written as  
\[
\eta = \gamma \left( \frac{\beta_2}{\beta_1(2^{\lambda - \log(\beta_1)}-1)} + 1 \right). \tag{64}
\]  
Consider $M=M(\eta)$. By Theorem 5 (or equivalently by Lemma 3 and Theorem 3), we have,
\[
E_{0}(M)>\frac{\gamma}{2}. \tag{65}
\]  
Hence, \( M \in S_\gamma \). For a choice of $\vartheta>0$ let $\lambda$ be chosen such that,
\[
{\lambda}=D(\mu_1 \| \hat{\mu}_0) - \frac{D(\mu_1 \| \mu_0)}{1 + \vartheta}. \tag{66}
\]
As \( \gamma \to \infty \), we have \( \epsilon_0 \to 0 \), hence \( n_0 \to \infty \), which yields  
\[
D(\mu_1 \| \hat{\mu}_0) \to D(\mu_1 \| \mu_0) \quad \text{a.s.} \ [\mu_0]. \tag{67}
\]  

Also, as \( n_0 \to \infty \), we have  
\[
\overline{E_1}(M\mid x_1^{n_0}) \to \overline{E_1}(M\mid x_1,x_2,x_3,...) \quad \text{a.s.} \ [\mu_0]. \tag{68}
\]    

From Theorem 5 we have,

\[\overline{E_1}(M \mid x_1^{n_0}) \sim \frac{\log \eta}{D(\mu_1 || \widehat{\mu_0}) - \lambda}. \tag{69}\]

Hence,
\[
\overline{E_1}(M \mid x_1,x_2,x_3,...) \sim \frac{\log \gamma}{D(\mu_1 || \widehat{\mu_0}) - \lambda} \quad \text{a.s.} \ [\mu_0]. \tag{70}
\]    
Substituting for $\lambda$,
\[
\overline{E_1}(M\mid x_1,x_2,x_3,...) \sim (1 + \vartheta)\frac{\log \gamma}{D(\mu_1 || {\mu_0})} \quad \text{a.s.} \ [\mu_0]. \tag{71}
\] 

From Theorem 4 (restatement of Theorem 4 of \cite{bansal1986}) we have,
\[\inf_{S \in S_\gamma} \overline{E_1}(S) \sim \frac{\log \gamma}{D(\mu_1 || \mu_0)}. \tag{72}\]

Hence, as \( n_0 \to \infty \) and because samples $x_1^{n_0}$ are generated using $\mu_0$, we get, 

\[
\overline{E_1}(M\mid x_1,x_2,x_3,...) \sim (1 + \vartheta) \inf_{S \in S_\gamma} \overline{E_1}(S) \quad \text{a.s.} \ [\mu_0]. \tag{73}
\] 

\section{Conclusion and Future Work}
The analysis holds for any stationary, ergodic, mutually independent finite-order Markov sources (pre-change distribution estimation, Lemma 2, Lemma 3 and equation (60) can be modified as per Markov source order). Hence, we have extended the results of 2022 for stationary and ergodic finite-order Markov sources under suitable conditions. One can aim to generalise these results further under Lai's framework \cite{lai1998} who has weakened the requirement of independence of observations before and after change to obtain an extension of Lorden's result. The performance of proposed test can also be studied under other alternate criteria introduced by Lai \cite{lai1998}. We have evaluated under Lorden's criteria as our work extends \cite{jacobbansal2008,jacob2010,malikbansal2022} all of which are studied under Lorden's criteria.

\end{document}